\documentclass[12pt]{article}
\usepackage{hyperref, graphicx}
 \usepackage{tikz}
 \usepackage{amsmath,amssymb,euscript,graphicx,amsfonts}
 \usepackage{textcomp}

\usepackage[english]{babel}

\def\proof{\par\noindent{\it Proof.\ \ }}
\def\qed{\ifmmode\square\else\nolinebreak\hfill
$\square$\fi\par\vskip12pt}

\newtheorem{theorem}{Theorem}[section]
\newtheorem{lemma}[theorem]{Lemma}%

\begin{document}
\baselineskip15pt

%\vskip 4in

\title{ On bipartite cages of excess $4$}
\author{
         Slobodan Filipovski\footnote{Supported in part by the Slovenian Research Agency (research program P1-0285 and Young Researchers Grant).} \\
         University of Primorska, Koper, Slovenia \\
         \texttt{slobodan.filipovski@famnit.upr.si}
 }

\maketitle

\begin{abstract}
The Moore bound $M(k,g)$ is a lower bound on the
order of $k$-regular graphs of girth $g$ (denoted $(k,g)$-graphs). The excess $e$ of a $(k,g)$-graph of order $n$ is the difference $ n-M(k,g) $. In this paper we consider the existence of $(k,g)$-bipartite graphs of excess $4$ via studying spectral properties of their adjacency matrices. We prove that the $(k,g)$-bipartite graphs of excess $4$ satisfy the equation $kJ=(A+kI)(H_{d-1}(A)+E)$, where $A$ denotes the adjacency matrix of the graph in question, $J$ the $n \times n$ all-ones matrix, $E$ the adjacency matrix of a union of
vertex-disjoint cycles, and $H_{d-1}(x)$ is the Dickson polynomial of the second kind with parameter $k-1$ and of degree $d-1$.
We observe that the eigenvalues other than $\pm k$ of these graphs are roots of the polynomials $H_{d-1}(x)+\lambda$, where $\lambda$ is an eigenvalue of $E$. Based on the irreducibility of $H_{d-1}(x)\pm2$ we give necessary conditions for the existence of these graphs.
 If $E$ is the adjacency matrix of a cycle of order $n$ we call the corresponding graphs \emph{graphs with cyclic excess}; if $E$ is the adjacency matrix of a disjoint union of two cycles we call the corresponding graphs \emph{graphs with bicyclic excess}. In this paper we prove the non-existence of $(k,g)$-graphs with cyclic excess $4$ if $k\geq6$ and $k \equiv1 \!\!  \pmod {3}$, $g=8, 12, 16$ or $k \equiv2 \!\!  \pmod {3}$, $g=8,$ and the non-existence of $(k,g)$-graphs with bicyclic excess $4$ if $k\geq7$ is odd number and $g=2d$ such that $d\geq4$ is even.

%\textbf{Keywords:} bipartite graphs, cage problem, girth

\end{abstract}


\section{Introduction}
\quad A $k$-regular graph of girth $g$ is called a $(k,g)$-{\em graph}. A $(k,g)$-\emph{cage} is a $(k,g)$-graph with the fewest possible number of vertices, among all $(k,g)$-graphs. The order of a $(k,g)$-cage is denoted by $n(k,g)$. The \textit{Cage Problem} calls for finding cages, and this problem was considered for the first time by Tutte \cite{K}. It is known that a $(k,g)$-graph exists for any combination of $k \geq2$ and $g\geq 3$,
\cite{erdossachs,sachs}. However, the orders $n(k,g)$ of $(k,g)$-cages have only been determined for very limited sets of parameters \cite{exojaj}.
A natural lower bound on the order of a $(k,g)$-graph is called the \textit{Moore bound}, and the form of the bound depends on the parity of $g$, i.e.,
\begin{equation} \label{moore1}
n(k,g) \ge  M(k,g) = \left\{
\begin{array}{lc} 1 + k + k(k-1) + ... + k(k-1)^{(g-3)/2},
                                                         & g \mbox{ odd, } \\
                              2 \left( 1 + (k-1) + ... + (k-1)^{(g-2)/2} \right),
                                                         & g \mbox{ even. }
\end{array}
\right.
\end{equation}

The graphs whose orders are equal to the Moore bound are called {\em Moore graphs}. They are known to exist
if $k=2$ and $ g \geq 3 $, $g=3$ and $ k \geq 2 $,
$g=4$ and $ k \geq 2 $, $g = 5$ and $k=2, 3, 7 $, or $g = 6, 8, 12$ and a generalized
$n$-gon of order $k-1$ exists \cite{banito,dam,exojaj}.
The existence of a $(57,5)$-Moore graph is an open question.\\
The {\em excess} $e$ of a $(k,g)$-graph is the difference between its order $n$ and the Moore
bound $M(k,g)$, i.e., $e = n - M(k,g)$.
 %For the majority
%of parameter pairs $(k,g)$ the exact values $n(k,g)$ are not known, and very few lower bounds on $n(k,g)$ exceeding the
%Moore bound exist.
%Brown \cite{brown} showed that $n(k,5)$ is never equal to $M(k,5)+1$.
%Eroh and Schwenk \cite{erosch} showed for $g=7$ the non-existence of
%$k$-regular graphs of girth $7$ and order $M(k,7)+1$.
%(in comparison, $n(3,7)=M(3,7)+2 $, the order of the McGee graph).
%%Regarding graphs with odd girth,
%Bannai and Ito \cite{Ban&Ito} have shown that no $k$-regular graphs of
%order $M(k,g)+1$ exist for any odd $g \geq 5$.
%Kov\' acs \cite{kovacs} has shown that no graphs of excess $2$, girth $5$, and odd
%degree $ k $ which is not of the form $ \ell^2 + \ell + 3 $ or $ \ell^2 + \ell - 1 $,
%where $ \ell $ is a positive integer, exist.
%Eroh and Schwenk \cite{erosch} proved that $n(k,5)$ is not equal to
%$M(k,5)+2$ for $5 \le k \le 11$. Most recent results concerning
%odd girth and excess $2$ are due to Garbe \cite{garbe}. He showed
%the non-existence of graphs of excess $2$ for parameters
%$(k,9)$, $ (k,13) $, $ (k,17) $, $ (k,21) $, $ (k, 25) $, and $ (k,29) $ for certain
%congruence classes of $k$.
%He also showed that there are no excess $2$ graphs in
%the families of
%$ (3,2s+1) $-graphs, $(7,2s+1)$-graphs, and $(9,2s+1)$-graphs,
%for certain congruence classes of $s$.
 Regarding graphs of even girth we will use the following three results:
\begin{theorem}[\cite{bigito}]\label{bigito}
Let $G$ be a $(k,g)$-cage of girth
$g = 2d \geq 6$ and excess $e$.
If $e \leq k-2$, then $e$ is even and $G$ is bipartite of diameter $d + 1$.
\end{theorem}
For the next theorem, let $D(k,2)$ denote the  incidence graph of a symmetric $(v,k,2)$-design.
\begin{theorem}[\cite{bigito}]\label{excess2}
Let $G$ be a a $(k,g)$-cage of girth $g = 2d \geq 6$ and excess $2$.
Then $g=6$, $G$ is a double-cover of $D(k,2)$, and $k$ is not congruent to $5$ or $7 \!\! \pmod{8}$.
\end{theorem}
\begin{theorem}[\cite{ourpaper}]\label{exclude}
Let $ k \geq 6 $, $ g = 2d > 6 $. No $(k,g)$-graphs of excess $4$ exist for parameters
$k,g$ satisfying at least one of the following conditions:
\begin{itemize}
\item[{\rm 1)}] $g=2p$, with $p\geq5$ a prime number, and $k \not \equiv 0, 1, 2 \!\! \pmod{p}$;
\item[{\rm 2)}] $ g=4\cdot3^{s}$ such that $ s\geq 4 $, and $k$ is divisible by $9$ but not by $ 3^{s-1}$;
\item[{\rm 3)}]  $ g=2p^{2}$ with $p\geq 5$ a prime number, and  $k \not \equiv 0, 1, 2 \!\! \pmod{p}$ and even;
\item[{\rm 4)}] $ g=4p$, with $p\geq 5 $ a prime number, and  $k \not \equiv 0, 1, 2,3,p-2 \!\! \pmod{p}$;
\item[{\rm 5)}] $g\equiv0\!\! \pmod {16}$, and $k\equiv3\!\! \pmod {g}.$
\end{itemize}
\end{theorem}

Motivated by the result in Theorem \ref{exclude}, which was obtained through counting cycles in a hypothetical graph with given parameters and excess $4$, in this paper we address the question of the existence of
 $(k,g)$-graphs of excess $4$ using spectral properties of their adjacency matrices. The question of the existence of $(k,g)$-graphs of excess $4$ is wide open, and prior to the publication of \cite{ourpaper}, no such results were known. The results contained in our paper further extend our understanding of the structure of the potential graphs of excess $4$.
Throughout, we assume that $k\geq6$, $g=2d\geq6$ and $G$ is a $(k,g)$-graph of excess $4$ and order $n$. Due to Biggs's result stated in Theorem \ref{bigito}, the restriction of the parameters $k$, $g$ given above allows us to conclude that $G$ is a bipartite graph with diameter $d+1$.
For each integer $i$ in the range $0\leq i\leq d+1$, we define the $n\times n$ matrix $A_{i}=A_{i}(G)$ as follows. The rows and columns of $A_{i}$ correspond to the vertices of $G$, and the entry in position $(u, v)$ is $1$ if the distance $d(u, v)$ between the vertices $u$ and $v$ is $i$, and zero otherwise. Clearly, $A_{0}=I, A_{1}=A$, the usual adjacency matrix of $G$. The last non-zero matrix is the matrix $A_{d+1}$ which we shall denote by $E$ and refer to it as the {\em excess matrix} i.e., $E$ is the adjacency matrix of the graph with the same vertex set $V$ as $G$  such that two vertices of $V$ are adjacent if and only if they have distance $d+1$. We will call this graph the \emph{excess graph} of $G$ and we will denote it by $G(E).$ If $J$ is the all-ones matrix, the sum of the \emph{$i$-distance matrices }$A_{i}, 0\leq i\leq d$, and the matrix $E$ yields $\sum_{i=0}^{d}A_{i}+E=J$. To apply the last identity we will use Lemma 4 from \cite{ourpaper}. Employing the methodology used by Bannai et al. in \cite{banito}, \cite{Ban&Ito}, later by Biggs et al. in \cite{bigito}, Delorme et al. in \cite{paper} and Garbe in \cite{garbe}, we will show that the eigenvalues of $G$ other than $\pm k$ are the roots of the polynomials $H_{d-1}(x)+\lambda$. Here, $H_{d-1}(x)$ is the Dickson polynomial of the second kind with parameter $k-1$ and degree $d-1$, and $\lambda$ is an eigenvalue of the excess matrix $E$. Furthermore, for odd $k\geq7$ and $d\geq4$, we prove that the polynomial $H_{d-1}(x)\pm2$ is irreducible over $\mathbb{Q}[x]$, which leads to necessary conditions for existence of $(k,g)$-graphs of excess $4$, Theorem \ref{a}.

   We say that a graph $G$ has a \emph{cyclic excess} if the excess graph $G(E)$ is a cycle of length $n$, and a graph $G$ has a \emph{bicyclic excess} if $G(E)$ is a disjoint union of two cycles.
In \cite{delorme} Delorme et al. considered graphs with cyclic defect and excess $2$, proving non-existence of infinitely many such graphs. The paper describes the cycle structure of the excess graphs of the known non-trivial graphs of excess $2$:
  \begin{itemize}
  \item[{1)}] the excess graph of the only $(3,5)$-graph of excess $2$ is a disjoint union of a $9$-cycle and a $3$-cycle or a disjoint union of an $8$-cycle and $4$-cycle;
  \item[{2)}] the excess graph of the unique $(4,5)$-graph of excess $2$ (the Robertson graph) is a disjoint union of a $3$-cycle, a $12$-cycle and a $4$-cycle;
  \item[{3)}] the excess graph of the unique $(3,7)$-graph of excess $2$ (the McGee graph) is a disjoint union of six $4$-cycles.
  \end{itemize}
 \quad We note that no $(k,g)$-graph of cyclic excess $2$ are known, while examples of graphs with bicyclic excess $2$ can be found among the $(3,5)$-graphs of excess $2$. Proving that the excess graphs of bipartite graphs of excess $4$ form a disjoint union of cycles, while also
 inspired by the results in \cite{delorme}, in Section $3$ we consider the existence of bipartite graphs of excess $4$ with cyclic and bicyclic excess $4$. Based on the irreducibility of $H_{d-1}(x)\pm2$ and $H_{d-1}(x)-1$ over $\mathbb{Q}[x]$, we prove the non-existence of infinitely many such graphs of girths at least $8$. %Considering the known bipartite graphs of excess $4$, let us mention the $(7,6)$-cage determined by O\textquotesingle Keefe and Wong \cite{unique}. It is known that this graph is vertex-transitive of order $n(7, 6) = 90=M(7,6)+4$.

%\section{The structure of graphs of even girth and excess $4$}
%In this section, we take on the case of
%$(k,g)$-graphs of degree $ k \geq 6 $, {\em even} girth $ g = 2m \geq 6 $, and excess $4$.
%All of these graphs are covered by Theorem~\ref{bigito} and are
%therefore bipartite and of diameter $ m+1 $. Thanks to these results,
%the structure of $G$ with respect to any edge $e = \{ u,v \}  \in E(G)$ can therefore be determined. Let $N_G(u,i)$ denote the
%$i$-th neighborhood of the vertex $u$, i.e., the set of vertices of $G$ whose distance from $u$ in $G$
%is equal to $i$. Since the girth of $G$ is assumed to be equal to $g$, the set of vertices of $G$

\section{Necessary conditions for the existence of graphs of even girth and excess $4$}

\quad Let $k\geq6$, $g=2d\geq6$, and let $G$ be a $(k,g)$-graph of excess $4.$ Then $G$ is bipartite of diameter $d+1.$
Let $N_G(u,i)$ denote the set of vertices of $G$ whose distance from $u$ in $G$
is equal to $i$, $1\leq i\leq d+1$. The
subgraph of $G$ induced by the set of vertices of $G$
whose distance from $u$ is at most $\frac{g-2}{2}$ and whose distance from $v$ is by one
larger than their distance from $u$ induces a tree of depth $\frac{g-2}{2}$ rooted at $u$ (we will  call
it ${\mathcal T}_u$). Also, the
subgraph of $G$ induced by the set of vertices of $G$
whose distance from $v$ is at most $\frac{g-2}{2}$ and whose distance from $u$ is by one
larger than their distance from $v$ induces a tree of depth $\frac{g-2}{2}$ rooted at $v$ (we will  call
it ${\mathcal T}_v$). Since $G$ is of girth $g$, the trees ${\mathcal T}_u$ and ${\mathcal T}_v$ are disjoint and contain no cycles.
Since each vertex of $G$ is of degree $k$, the order of ${\mathcal T}_u\bigcup{\mathcal T}_v$ is equal to
$ 2(1 + (k-1) + (k-1)^2 + \ldots +(k-1)^{\frac{g-2}{2}}) $.
%The degrees of $u$ or $v$ in their respective trees are equal to $(k-1)$, the
%degrees of all the non-leave vertices of these trees are equal to $k$, and all the leaves of these
%trees are of distance $\frac{g-2}{2}$ from their respective roots. As for the order of these subtrees,
%they are both of order $ 1 + (k-1) + (k-1)^2 + \ldots +(k-1)^{\frac{g-2}{2}} $, with $(k-1)^i$ vertices
%of distance $i$ from $u$ (or $v$).
We will call the union of the trees ${\mathcal T}_u, {\mathcal T}_v $ with the edge $f$ {\em Moore tree of $G$ rooted at $f$}; it is the subtree of $G$ that
is the basis of the Moore bound for even $g$. The graph $G$
must contain $4$ additional vertices $ w_1,w_2,w_3,w_4 $ which do not belong to either $ {\mathcal T}_u $
or $ {\mathcal T}_v $, and whose distance from both $u$ and $v$ is greater than $\frac{g-2}{2}$.
We will call these vertices {\em the excess vertices with respect to $f$} and denote this set
$ X_f = \{ w_1,w_2,w_3,w_4 \} $; we call the edges not contained in the Moore tree
of $G$ {\em horizontal edges}.

The following lemma restricts the possible ways in which the
four excess vertices are attached to the Moore tree.

\begin{lemma}[\cite{ourpaper}] \label{lemma4}
Let $ k \geq 6 $, $ g = 2d \geq 6 $. Let $G$ be a $(k,g)$-graph
of excess $4$, $u,v$ be two adjacent vertices in $G$, and $ X_f = \{ w_1,w_2,w_3,
w_4 \} $
be the four excess vertices with respect to the edge $ f = \{ u,v \} $. The
induced subgraph $ G[w_1,w_2,w_3,w_4] $ is isomorphic to $ 2K_2$
(two disjoint copies of $K_2$) or ${\mathcal P}_3$ (a path of length $3$).
\end{lemma}
Next, let us define the following polynomials:
\begin{center}
$F_{0}(x)=1, F_{1}(x)=x, F_{2}(x)=x^{2}-k;$
\end{center}
\begin{center}
$G_{0}(x)=1, G_{1}(x)=x+1;$
\end{center}
\begin{center}
$H_{-2}(x)=-\frac{1}{k-1}, H_{-1}(x)=0, H_{0}(x)=1, H_{1}(x)=x;$
\end{center}

\begin{equation}\label{2}
P_{i+1}(x)=xP_{i}(x)-(k-1)P_{i-1}(x) \mbox{ for }
\left\{
\begin{array}{lc} i\geq2, & \mbox{ if } P=F, \\
i\geq1, & \mbox{ if } P=G, \\
i\geq1, & \mbox{ if } P=H.
\end{array}
\right.
\end{equation}

In \cite{singleton}, Singleton gives many relationships between these polynomials. We will use two of them. Given any $i\geq0$,
\begin{eqnarray}\label{eq3}
G_{i}(x)=\sum_{j=0}^{i}F_{j}(x)
\end{eqnarray}
\begin{eqnarray}\label{eq4}
G_{i+1}(x)+(k-1)G_{i}(x)=(x+k)H_{i}(x).
\end{eqnarray}

The above defined polynomials have a close connection to the properties of a graph $G$. Namely, for $ t < g $ the element $(F_{t}(A))_{x,y}$ counts the number of paths of length $t$ joining vertices $x$ and $y$ of $G$. It follows from (\ref{eq3}) that $G_{t}(A)$ counts the number of paths of length at most $t$ joining pairs of vertices in $G$. All of the preceding claims can be found in \cite{paper}.
 The next lemma is based on the structure of $G$ described in Lemma \ref{lemma4}:
\begin{lemma}\label{vazna} Let $k\geq6$, $g=2d\geq6$ and let $G$ be a $(k,g)$-graph of excess $4$. If $A$ is the adjacency matrix of $G$ and $E$ is the excess matrix of $G$, then
$$F_{d}(A)=kA_{d}-AE.$$
\end{lemma}
\noindent

\proof Let $f=\{u,v\}$ be a base edge of the Moore tree and let $f_{1}=\{w_{1}, w_{2}\}, f_{2}=\{w_{3}, w_{4}\}$ be the edges of the subgraph induced by $X_f$. Also, let us assume that $d(u,w_{1})=d(u,w_{3})=d$ and
$d(u,w_{3})=d(u,w_{4})=d+1.$ We consider the case when $ G[w_1,w_2,w_3,w_4]$ is isomorphic to $2K_{2}$ in which case the excess vertices do not share common neighbour. The other cases
when $G[w_1,w_2,w_3,w_4]$ is isomorphic to $2K_{2}$ and the excess vertices share common neighbour or the subgraph induced by the excess vertices contains $ {\mathcal P}_3 $ are analogous. Since there are $k-1$ paths of length $d$ from $u$ to $w_{1}$ and $w_{3}$, by the definition of $F_{i}(x)$ we have $(F_{d}(A))_{u,w_{1}}=(F_{d}(A))_{u,w_{3}}=k-1.$ Considering the vertices of distance $d$ from $u$, there are also the $(k-1)^{d-1}$ leaves of the subtree rooted at $v$. For $2(k-1)$ of these vertices there exists $k-1$ paths of length $d$ from $u$ to them. Namely, they are the vertices adjacent to $w_{2}$ or $w_{4}.$ For all the other leaves, there are $k$ paths between. Thus, $(F_{d}(A))_{u,s}=0$ if $d(u,s)\neq d$, $(F_{d}(A))_{u,s}=k$ if $s$ is a leaf of a branch rooted at $v$ and not adjacent to $w_{2}$ and $w_{4}$, and $(F_{d}(A))_{u,s}=k-1$ if $s$ is $w_{1}, w_{3}$ or a leaf of a branch rooted at $v$ and adjacent to $w_{3}$ or $w_{4}$. This yields for the matrix $kA_d$ that $(kA_{d})_{u,s}=k$ if $d(u,s)=d$ and $(kA_{d})_{u,s}=0$ if $d(u,s)\neq d.$
Now, let $s$ be a vertex of $G$ such that $d(u,s)=d$ and $s$ is adjacent to $w_{2}$ or $w_{4}.$ If $s=w_{1}$ or $s=w_{3}$ then it is easy to see that $(AE)_{u,s}=1.$ On the other hand, since $s$ is adjacent to the subtree rooted at $u$ through $k-2$ different horizontal edges, it follows that between the $k-1$ branches of the subtree rooted at $u$ there exists one sub-branch that is not adjacent to $s$ though a horizontal edge. Let $s_{1}$ be the root of that sub-branch. Then, $d(s,s_{1})=d+1, d(u,s_{1})=1$, which implies $(A)_{u,s_{1}}=1$ and $(E)_{s_{1}, s}=1.$
Let $s_{2}$ be the other vertex of distance $d+1$ from $s$. Because all neighbours of $u$, except $s_{1}$, are of distance smaller than $d+1$ of $s$, we have $(A)_{u,s_{2}}=0$ and $(E)_{s_{2}, s}=1.$ Thus $(AE)_{u,s}=1.$
 If $s$ is a vertex of $G$ such that $d(u,s)=d$ and $s$ is not adjacent to $w_{2}$ or $w_{4}$ then the distance between $s$ and the neighbours of $u$ is $d-1$. In this case, $(AE)_{u,s}=0.$
If $d(u,s)\neq d$ then the distance between $s$ and the neighbours of $u$ is different from $d+1$, and therefore $(AE)_{u,s}=0. $
The required identity follows from summing up the above conclusions.
\hfill {\em q.e.d.}

\begin{lemma} Let $k\geq6$, $g=2d\geq6$ and let $G$ be a $(k,g)$-graph of excess $4$. If $A$ is the adjacency matrix of $G$, $E$ is the excess matrix of $G$ and $J$ is the all-ones matrix, then
$$kJ=(A+kI)(H_{d-1}(A)+E).$$
\end{lemma}

\noindent

\proof By the definition of the polynomials $G_{i}(x)$ and using the fact that $G$ has diameter $d+1$ we conclude $J=G_{d-1}(A)+A_{d}+E.$
The relation (\ref{eq3}), setting $i=d$, asserts $G_{d}(A)=G_{d-1}(A)+F_{d}(A)$.
Substituting this identity in (\ref{eq4}), where we fix $i=d-1$, we get $ kG_{d-1}(A)+F_{d}(A)=(A+kI)H_{d-1}(A).$
Due to Lemma \ref{vazna} the last identity is equivalent to $ kG_{d-1}(A)+kA_{d}+kE=(A+kI)(H_{d-1}(A)+E).$
From $kJ=kG_{d-1}(A)+kA_{d}+kE$ follows $kJ=(A+kI)(H_{d-1}(A)+E).$
\hfill {\em q.e.d.}
\medskip

The next theorem gives a relationship between the eigenvalues of the matrices $A$ and $E$ (this result is an analogue of Theorem 3.1 in \cite{paper}):
\begin{theorem} \label{eigenvalues}
If $\mu(\neq\pm k)$ is an eigenvalue of $A$, then
$$H_{d-1}(\mu)=-\lambda,$$
where $\lambda$ is an eigenvalue of $E$.
\end{theorem}
\noindent

\proof Let us suppose that $\mu$ is an eigenvalue of $A$. Since $G$ is a $k$-regular graph, the all-ones matrix $J$ is a polynomial in $A$.
  This implies that any eigenvector of $A$ is also an eigenvector of $J$. From $kJ=(A+kI)(H_{d-1}(A)+E)$ and since $H_{d-1}(A)$ is also a polynomial in $A$, we have that $E$ is a polynomial in $A$, and consequently, every eigenvector of $A$ is an eigenvector of $E$. Therefore, the eigenvalues of $kJ$ are of the form $(\mu+k)(H_{d-1}(\mu)+\lambda)$. As is well-known, the eigenvalues of $kJ$ are $kn$ (with multiplicity $1$) and $0$ (with multiplicity $n-1$).
  The eigenvalue $kn$ corresponds to $\mu=k$, and so all the remaining eigenvalues, except for $-k$, satisfy the above equation.
\hfill {\em q.e.d.}

Since the eigenvalues of a disjoint union of cycles are known, we are now in a position to determine the spectrum of $A$:

\begin{lemma}\label{spektar} Let $k\geq6$, $g=2d\geq6$ and let $G$ be a $(k,g)$-graph of excess $4$. If $A$ is the adjacency matrix of $G$ and $E$ is the excess matrix of $G$, then:
\begin{itemize}
\item[{1)}] The matrix $E$ is the adjacency matrix of a graph $G(E)$, consisting of a disjoint union of $c$ cycles $C_{i}$ of length $l_{i}$ with $1\leq i\leq c$. Moreover, if $d$ is odd and $V_{1}$ and $V_{2}$ are the two partition sets of the bipartite graph $G$, then every cycle in $G(E)$ is completely contained either in $V_{1}$ or $V_{2}$.

\item[{2)}] The spectrum of $A$ consists of:
\begin{itemize}
\item [{2.1)}] $\pm k$, $c-2$ many solutions of $H_{d-1}(x)=-2$, and one solution of each equation $H_{d-1}(x)=-2\cos(\frac{2\pi j}{l_{i}}), j=1,..., l_{i}-1; 1\leq i\leq c$, for $d$ odd;
    \item [{2.2)}]  $\pm k$, $c-1$ many solutions of $H_{d-1}(x)=-2$, and one solution of each equation (except one) $H_{d-1}(x)=-2\cos(\frac{2\pi j}{l_{i}}), j=1,..., l_{i}-1; 1\leq i\leq c$, for $d$ even;

        \end{itemize}
\end{itemize}
\end{lemma}
\noindent

\proof $1)$ Our proof is analogous to that of Kov\' acs for girth $5$, \cite{kovacs}, and Garbe's proof for odd girth $g=2k+1>5$, \cite{garbe}.
Let $f=\{u, v\}$ be a base edge of a bipartite Moore tree of $G$. Lemma \ref{lemma4} asserts that there exist exactly two vertices of $G$ on distance $d+1$ from $u$. Namely, they are the excess vertices adjacent to the leaves of the subtree rooted at $v$. The excess matrix $E$ is the adjacency matrix for the graph $G(E)$ with same vertex set $V$ as $G$ such that two vertices of $G(E)$ are adjacent if and only if they are of distance $d+1$. Because for each vertex $u\in V(G)$ there are exactly two vertices on distance $d+1$ from $u$, every component of $G(E)$ is a cycle. Let $c$ be the number of these cycles and let $l_{i}, i=1,..,c$, be the lengths of these cycles ordered in an arbitrary manner. Moreover, if $d$ is an odd number, any two vertices of $G$ with distance $d+1$ lie in the same partite set. Therefore any connected component of $G(E)$ is entirely contained either in $V_{1}$ or $V_{2}$.\\
$2)$ The eigenvalues of an $n$-cycle are known and are equal to $2\cos(\frac{2\pi j}{n}), (j=0,...,n-1)$.
Therefore the eigenvalues of $G(E)$ are $2\cos(\frac{2\pi j}{l_{i}}), j=0, 1,..., l_{i}-1; 1\leq i\leq c$, \cite{garbe}.
Since $G$ is a $k$-regular bipartite graph, it has (among others)
the eigenvalues $k$ and $-k.$ Let $V_{1}$ and $V_{2}$ be the partition sets of $G$. Hence the eigenvector of $A$
corresponding to $k$ consist of the all-ones vector $j$, and the eigenvector corresponding to $-k$
is the vector $j'$ with values $1$ on $V_1$ and values $-1$ on $V_2$. If $d$ is an odd number then two vertices of $G(E)$ are adjacent if and only if they are in the same partite set. Therefore $E\cdot j'=2j'$, which implies that from the set of $c$ solutions on $H_{d-1}(x)=-2$ we need to subtract two multiplicities for the eigenvalues $k$ and $-k$ of $G$. If $d$ is an even number then two vertices of $G(E)$ are adjacent if and only if they are in different partite sets. Thus $E\cdot j'=-2j'$.  In this case, from the set of $c$ solutions on $H_{d-1}(x)=-2$ we need to subtract one multiplicity for the eigenvalue $k$ and from the set of all solutions on $H_{d-1}(x)=2$ we need to subtract one multiplicity for the eigenvalue $-k$. \hfill {\em q.e.d.}

\begin{lemma} \label{important} Let $k\geq6$, $g=2d\geq6$ and let $G$ be a $(k,g)$-graph of excess $4$. Furthermore, let $c$ be the number of cycles of $G(E)$ and $c_{2}$ be the number of cycles of even length. Then:
\begin{itemize}
\item[{1)}] If $H_{d-1}(x)-2$ is irreducible over $\mathbb{Q}[x]$ then $d-1$ divides $c-1$ or $c-2$;
\item[{2)}] If $H_{d-1}(x)+2$ is irreducible over $\mathbb{Q}[x]$, then $d-1$ divides $c_{2}-1$ or $c_{2}.$
\end{itemize}
\end{lemma}
\noindent
\proof
 $1)$ Combining Theorem \ref{eigenvalues} and part 2) from Lemma \ref{spektar} we obtain that $H_{d-1}(x)-2$ is an irreducible factor of the characteristic polynomial of $A$. Realizing that the roots of an irreducible factor of a characteristic polynomial of given rational symmetric matrix have the same multiplicities, \cite{kovacs}, from $2)$ of Lemma \ref{spektar} we have:\\
     If $d$ is an even number then the $d-1$ roots of $H_{d-1}(x)-2$ have multiplicity $\frac{c-1}{d-1}$, which has to be a positive integer.
      If $d$ is odd then the $d-1$ roots have multiplicity $\frac{c-2}{d-1}$.\\
      $2)$ Part 2) follows along the same lines as part 1).
      \hfill {\em q.e.d.}\\

      \medskip

We can base the testing of irreducibility of $H_{d-1}(x)\pm2$ on the well-known Eisenstein's criterion that asserts for a polynomial $f(x)=\sum_{i=0}^{n}a_{i}x^{i}\in \mathbb{Z}[x]$ and a prime $p$ that divides $a_{i}$ for all $0\leq i<n$, does not divide $a_{n}$ and $p^{2}$  does not divide  $a_{0}$ Now we are ready for the main result in this section:
\begin{theorem}\label{a} Let $k\geq7$ be an odd number and let $g=2d\geq8$. Let $c$ be the number of cycles of $G(E)$ and $c_{2}$ be the number of cycles with even length. If there exists a $(k,g)$-graph of excess $4$ then

\begin{itemize}
   \item[{1)}] if $d$ is an odd number then $d-1$ divides $c-2$ and $c_{2};$
      \item[{2)}] if $d$ is an even number then $d-1$ divides $c-1$ and $c_{2}-1.$
      \end{itemize}
 \end{theorem}
\noindent
\proof According to Lemma \ref{important}, it is enough to prove that the polynomials $H_{d-1}(x)-2$ and $ H_{d-1}(x)+2$ are irreducible.
We will prove using induction on $d\geq4$ that $H_{d-1}(x)=x^{d-1}+(k-1)\cdot P_{d-3}(x)$,
   where $P_{d-3}(x)$ is an integer polynomial of degree $d-3$. For $d=4$ we calculate $H_{3}(x)=x^{3}-2(k-1)x.$ Let us suppose that the above formula holds for $H_{d-2}(x)$ and $H_{d-3}(x)$. That yields
\begin{center}
$H_{d-1}(x)=x(x^{d-2}+(k-1)\cdot P_{d-4}(x))-(k-1)(x^{d-3}+(k-1)\cdot P_{d-5}(x))$=
$=x^{d-1}+(k-1)\cdot P_{d-3}(x).$
\end{center}
Therefore $H_{d-1}(x)\pm2=x^{d-1}+(k-1)\cdot P_{d-3}(x)\pm2.$
By the inductional hypothesis, follows that for an odd $d$ occurs $H_{d-1}(0)=(-1)^\frac{d-1}{2}\cdot (k-1)^\frac{d-1}{2}$
and $H_{d-1}(0)=0$ for an even $d$.
Hence for an odd $d\geq5$ the absolute value $(-1)^\frac{d-1}{2}\cdot (k-1)^\frac{d-1}{2}\pm 2$ is not divisible by $2^{2}$, and clearly
for an even $d\geq4$, $\pm2$ is not divisible by $2^{2}.$
Since $k-1$ is even, it follows that every coefficient on $H_{d-1}(x)\pm2$ except for the coefficient $1$ of $x^{d-1}$ is divisible by $2$. Thus, the conditions of the Eisenstein's criterion are satisfied, and $H_{d-1}(x)\pm2$ is irreducible.
\hfill {\em q.e.d.}
\section{The non-existence of bipartite graphs of cyclic or bicyclic excess}

\quad In this section we still deal with the family of graphs considered as in Section $2$. Again, let $k\geq6, g=2d\geq6$ and let $G$ be a $(k,g)$-graph of excess $4$ and order $n$. Clearly $n$ is even number.
We have already proved that the excess graph $G(E)$ consists of a disjoint union of $c$ cycles $C_{i}, 1\leq i\leq c$.
If $c=1$ and $G(E)$ consists of an $n$-cycle, $G$ is of cyclic excess $4$, and if $c=2$ and $G(E)$ consists of a disjoint union of two cycles, $G$ is of bicyclic excess $4$. These are the graphs we study in this section.
Note that there are no graphs $G$ with cyclic excess $4$ if $d$ is an odd number; in this case we showed that each cycle of $G(E)$ is completely contained either in $V_{1}$ or $V_{2}.$

Let $d$ be an even number and let $L_{n}$ be an $n$-cycle formed by the vertices of $G(E)$. If $A^{'}$ is the adjacency matrix of $L_{n}$, its characteristic polynomial $\chi(L_{n}, x)$ satisfies $\chi(L_{n},x)=(x-2)(x+2)(R_{n}(x))^{2}$,
where $R_{n}$ is a monic polynomial of degree $\frac{n}{2}-1$.
Consider the factorization $x^{n}-1=\prod_{l|n}\Phi_{l}(x),$ where $\Phi_{l}(x)$ denotes the $l$-\emph{th cyclotomic polynomial}. In the following paragraph, we summarize the properties of cyclotomic polynomials as listed in \cite{delorme}.\\
The cyclotomic polynomial $\Phi_{l}(x)$ has integral
coefficients, it is irreducible over $\mathbb{Q}[x]$, and it is \emph{self-reciprocal} ($x^{\phi(l)}\Phi_{l}(1/x)=\Phi_{l}(x)).$
From the irreducibility and the self-reciprocity of $\Phi_{l}(x)$ follows that the degree of $\Phi_{l}(x)$ is even for $l\geq2.$\\
Thus we obtain the following factorization of $R_{n}(x): R_{n}(x)=\prod_{3\leq l|n}f_{l}(x),$ where $f_{l}$ is an integer polynomial of degree $\frac{\phi(l)}{2}$ satisfying $x^{\phi(l)/2}f_{l}(x+1/x)=\Phi_{l}(x).$
Also, $f_{l}$ is irreducible over $\mathbb{Q}[x]$ and $f_{3}(x)=x+1, f_{4}(x)=x, f_{5}(x)=x^{2}+x-1, f_{6}(x)=x-1.$
Substituting $y=-H_{d-1}(x)$ into $\frac{\chi(L_{n}, y)}{(y-2)}$, we obtain a polynomial $F(x)$ of degree $(n-1)(d-1)$ which satisfies $F(A)u=0$ for each eigenvector $u$ of $A$ orthogonal to the all $-1$ vector. Setting $F_{l,k,d-1}(x)=f_{l}(-H_{d-1}(x))$ yields
$$F(x)=(-H_{d-1}(x)+2)\prod_{3\leq l|n}(F_{l,k,d-1}(x))^{2}.$$

\begin{lemma} \label{bitna} Let $g=2d>6$ and $l\geq3$ be a divisor of $n$. If there is a $(k,g)$-graph with cyclic excess $4$ and order $n$, then $F_{l,k,d-1}(x)$ must be reducible over $Q[x]$.
\end{lemma}

\noindent
\proof The degree of $F_{l,k,d-1}(x)$ is equal to $(d-1)\cdot\frac{\phi(l)}{2}$. If $F_{l,k,d-1}(x)$ is irreducible over $\mathbb{Q}[x],$ then all its roots must be eigenvalues of $A$. Employing Observation 3.1. from \cite{delorme}, we conclude that there are at most $\phi(l)$ roots of $F_{l,k, d-1}(x)$ that are eigenvalues of $A$.
Thus $(d-1)\cdot\frac{\phi(l)}{2}$=$\phi(l)$ i.e., $d=3$. This contradicts the assumption that $2d > 6$.
\hfill {\em q.e.d.}\\

Note that $deg(F_{l,k,d-1}(x))=d-1$ if and only if $\phi(l)=2,$ i.e., if and only if $l\in\{3, 4, 6\}.$

\begin{lemma} \label{glavna} Let $k\geq6, g=2d>6$, and let $n$ be the order of a $(k,g)$-graph with cyclic excess $4.$ Then

\begin{itemize}
\item[{\rm 1)}]if $n\equiv 0 \!\!  \pmod {3}$, then $H_{d-1}(x)-1$ must be reducible over $\mathbb{Q}[x]$;
 \item[{\rm 2)}]if $n\equiv 0 \!\!  \pmod {4}$, then $H_{d-1}(x)$ must be reducible over $\mathbb{Q}[x];$
\item[{\rm 3)}]   if $n\equiv 0 \!\!  \pmod {6}$, then $H_{d-1}(x)+1$ must be reducible over $\mathbb{Q}[x].$
\end{itemize}
\end{lemma}
\noindent
\proof Follows directly from Lemma \ref{bitna}, with the additional assumption $f_{3}(x)=x+1, f_{4}(x)=x$ and $f_{6}(x)=x-1$.
\hfill {\em q.e.d.}

If $n\equiv0 \!\!  \pmod {4}$, then using the formula for the order of $G$, $d-1$ must be odd.
On the other hand, since $H_{1}(x)=x, H_{3}(x)=x^{3}-2(k-1)x$ and $H_{d-1}(x)=xH_{d-2}(x)-(k-1)H_{d-3}(x)$  we see that if $d-1$ is an odd number then $x$ divides $H_{d-1}(x)$, which implies that $H_{d-1}(x)$ is reducible.
Therefore the second condition from Lemma \ref{glavna} is satisfied.\\
\quad The irreducibility of the polynomials $H_{d-1}(x)-1$ over $\mathbb{Q}[x]$ is examined in \cite{paper}, where it is analytically proven that these polynomials are irreducible for $d\in\{4,6,8\}$ and the paper contains a conjecture that $d\geq10$, $H_{d-1}(x)-1$ is irreducible.
From the irreducibility of $H_{d-1}(x)-1$ we obtain the main non-existence result of our paper.

\begin{theorem} If $k$ and $g$ satisfy one of the following conditions, there exist no $(k,g)$-graphs of cyclic excess $4$:
\begin{itemize}
\item[{\rm 1)}] $k\equiv1, 2 \!\! \pmod {3}$ and $g=8$;
\item[{\rm 2)}] $k\equiv1 \!\!  \pmod {3}$ and $g=12$;
\item[{\rm 3)}] $k\equiv1 \!\!  \pmod {3}$ and $g=16$.
\end{itemize}
\end{theorem}

\noindent
\proof Because the order of the graphs is equal to $4+2 \left( 1 + (k-1) + ... + (k-1)^{(g-2)/2}\right)$ we conclude $n\equiv 0 \!\! \pmod {3}$.
 Since the polynomial $H_{d-1}(x)-1$ is known to be irreducible for $d\in\{4, 6, 8\}$, we get contradiction to 1) from Lemma \ref{glavna}.
\hfill {\em q.e.d.}

\emph{Remark: }Since $d$ is an even number, Theorem \ref{a} asserts that $d-1$ divides $c-1$ and $c_{2}-1.$
This claim is satisfied because $c=c_{2}=1.$

Next, let us consider graphs of bicyclic excess $4$. In this case, we can assume an arbitrary (even or odd) $d$, as this case does not depend of the parity of $d$. So, let $G(E)$ be a graph consisting of a disjoint union of two cycles $C_{1}$ and $C_{2}$. If $d$ is an odd number, then the vertex sets of the cycles $C_{1}$ and $C_{2}$ correspond to the partite sets $V_{1}$ and $V_{2}$, respectively.\\
If $n\equiv0 \!\!  \pmod {4}$, $d$ is an even, each edge of $C(E)$ has endpoints in $V_{1}$ and $V_{2}$, and therefore each of the cycles has even length, $c_{2}=2.$ Furthermore, $k-1$ must be odd. Unfortunately, this will not help us in excluding any family of pairs $(k,g)$ for which $G$ does not exist. In fact, for an odd $d-1$ and an odd $k-1$ we cannot conclude irreducibility of $H_{d-1}(x)+2,$ thus, we cannot employ Lemma \ref{important}.\\
If $n\equiv2 \!\! \pmod {4}$ and $d$ is odd, then the lengths of $C_{1}$ and $C_{2}$ are equal to $\frac{n}{2}$ (clearly $n=2s+1$ is odd). Therefore $c_{2}=0$
and clearly $d-1$ divides $c-2$ and $c_{2}.$

The main result about the non-existence of graphs $G$ with bicyclic excess $4$ is given in the following theorem:
\begin{theorem}\label{b}  If $k\geq7$ is an odd and $g=2d\geq8$, where $d$ is an even integer, then there exist no $(k,g)$-graphs with bicyclic excess $4$.
\end{theorem}

\noindent
\proof We have $c=2.$ Theorem \ref{a} implies that $d-1$ divides $c-1$; a contradiction.
\hfill {\em q.e.d.}

\end{document}